\documentclass[english,12pt]{smfart}
\usepackage{amssymb}
\usepackage{amsfonts}
\usepackage{amscd}

\swapnumbers

\def\ac{{\rm ac}}
\def\pr{{\rm pr}}

\def\Var{{\rm Var}}

\def\GM{{\mathbf G}_m}

\def\AA{{\mathbf A}}

\def\LL{{\mathbf L}}

\def\NN{{\mathbf N}}

\def\RR{{\mathbf R}}

\def\ZZ{{\mathbf Z}}

\def\cA{{\mathcal A}}

\def\cC{{\mathcal C}}

\def\cI{{\mathcal I}}

\def\cL{{\mathcal L}}
\def\cM{{\mathcal M}}

\def\cO{{\mathcal O}}

\def\cS{{\mathcal S}}

\def\cX{{\mathcal X}}
\def\cY{{\mathcal Y}}

\mathchardef\alphag="7C0B \mathchardef\betag="7C0C
\mathchardef\gammag="7C0D \mathchardef\deltag="7C0E
\mathchardef\varepsilong="7C22 \mathchardef\varphig="7C27
\mathchardef\psig="7C20 \mathchardef\zetag="7C10
\mathchardef\epsilong="7C0F \mathchardef\rhog="7C1A
\mathchardef\taug="7C1C \mathchardef\upsilong="7C1D
\mathchardef\iotag="7C13 \mathchardef\thetag="7C12
\mathchardef\pig="7C19 \mathchardef\sigmag="7C1B
\mathchardef\etag="7C11 \mathchardef\omegag="7C21
\mathchardef\kappag="7C14 \mathchardef\lambdag="7C15
\mathchardef\mug="7C16 \mathchardef\xig="7C18
\mathchardef\chig="7C1F \mathchardef\nug="7C17
\mathchardef\varthetag="7C23 \mathchardef\varpig="7C24
\mathchardef\varrhog="7C25 \mathchardef\varsigmag="7C26
\mathchardef\Omegag="7C0A \mathchardef\Thetag="7C02
\mathchardef\Sigmag="7C06 \mathchardef\Deltag="7C01
\mathchardef\Phig="7C08 \mathchardef\Gammag="7C00
\mathchardef\Psig="7C09 \mathchardef\Lambdag="7C03
\mathchardef\Xig="7C04 \mathchardef\Pig="7C05
\mathchardef\Upsilong="7C07

\newtheorem{theorem}[subsection]{Theorem}
\newtheorem{lem}[subsection]{Lemma}

\theoremstyle{definition}

\newtheorem{example}[subsection]{Example}

\newtheorem{def-prop}[subsubsection]{Definition-Proposition}

\theoremstyle{remark}
\newtheorem{remark}[subsection]{Remark}

\theoremstyle{plain}

\numberwithin{equation}{subsection}

\def\boxit#1#2{\setbox1=\hbox{\kern#1{#2}\kern#1}%
\dimen1=\ht1 \advance\dimen1 by #1 \dimen2=\dp1
\advance\dimen2 by #1
\setbox1=\hbox{\vrule height\dimen1 depth\dimen2\box1\vrule}%
\setbox1=\vbox{\hrule\box1\hrule}%
\advance\dimen1 by .4pt \ht1=\dimen1 \advance\dimen2 by
.4pt \dp1=\dimen2 \box1\relax}

\newcommand{\sur}[2]{\genfrac{}{}{0pt}{}{#1}{#2}}
\renewcommand{\theequation}{\thesubsection.\arabic{equation}}

\def\AA{{\mathbf A}}

\def\LL{{\mathbf L}}

\def\NN{{\mathbf N}}

\def\RR{{\mathbf R}}

\def\ZZ{{\mathbf Z}}

\def\cA{{\mathcal A}}

\def\cC{{\mathcal C}}

\def\cI{{\mathcal I}}

\def\cL{{\mathcal L}}
\def\cM{{\mathcal M}}

\def\cO{{\mathcal O}}

\def\cS{{\mathcal S}}

\def\cX{{\mathcal X}}
\def\cY{{\mathcal Y}}

\mathchardef\alphag="7C0B \mathchardef\betag="7C0C
\mathchardef\gammag="7C0D \mathchardef\deltag="7C0E
\mathchardef\varepsilong="7C22 \mathchardef\varphig="7C27
\mathchardef\psig="7C20 \mathchardef\zetag="7C10
\mathchardef\epsilong="7C0F \mathchardef\rhog="7C1A
\mathchardef\taug="7C1C \mathchardef\upsilong="7C1D
\mathchardef\iotag="7C13 \mathchardef\thetag="7C12
\mathchardef\pig="7C19 \mathchardef\sigmag="7C1B
\mathchardef\etag="7C11 \mathchardef\omegag="7C21
\mathchardef\kappag="7C14 \mathchardef\lambdag="7C15
\mathchardef\mug="7C16 \mathchardef\xig="7C18
\mathchardef\chig="7C1F \mathchardef\nug="7C17
\mathchardef\varthetag="7C23 \mathchardef\varpig="7C24
\mathchardef\varrhog="7C25 \mathchardef\varsigmag="7C26
\mathchardef\Omegag="7C0A \mathchardef\Thetag="7C02
\mathchardef\Sigmag="7C06 \mathchardef\Deltag="7C01
\mathchardef\Phig="7C08 \mathchardef\Gammag="7C00
\mathchardef\Psig="7C09 \mathchardef\Lambdag="7C03
\mathchardef\Xig="7C04 \mathchardef\Pig="7C05
\mathchardef\Upsilong="7C07

\DeclareMathOperator*{\Spec}{Spec}
\def\ord{{\rm ord}}

\begin{document}

\title[Nearby cycles
and non-degenerate polynomials]{Nearby cycles
and composition with a non-degenerate polynomial}

\author{Gil Guibert}

\address{39 quai du Halage,
94000 Cr\'eteil, France}

\email{guibert9@wanadoo.fr}

\author{Fran\c cois Loeser}

\address{{\'E}cole Normale Sup{\'e}rieure,
D{\'e}partement de math{\'e}matiques et applications, 45
rue d'Ulm, 75230 Paris Cedex 05, France (UMR 8553 du CNRS)}
\email{Francois.Loeser@ens.fr}
\urladdr{http://www.dma.ens.fr/\~{}loeser/}

\author{Michel Merle}

\address{Laboratoire J.-A. Dieudonn\'e,
Universit\'e de Nice - Sophia Antipolis, Parc Valrose,
06108 Nice Cedex 02, France (UMR 6621 du CNRS)}
\email{Michel.Merle@unice.fr}
\urladdr{http://www-math.unice.fr/membres/merle.html}
%\dedicatory{Preliminary version, October 30, 2001}

%\begin{abstract}We express the Lefschetz number of iterates of the
%monodromy of a function on a smooth complex algebraic variety
%in terms of the Euler characteristic of a space of truncated arcs.
%\end{abstract}
\maketitle

\section{Introduction}Let $X_j$ be smooth varieties
over a field $k$ of characteristic zero, for $1\leq j\leq p$.
Consider a family $\mathbf{f}$ of $p$  functions $f_j : X_j
\rightarrow \AA^1_k$. We shall denote also by
$f_j$ the function on the product $X=\prod_j X_j$
obtained by composition with the projection. We denote by 
$X_0(\mathbf{f})$ the set of  common zeroes in $X$ of the functions $f_j$.
Let $P\in k[y_1,\ldots,y_p]$ be a polynomial, which we assume to be
non-degenerate with respect to its Newton polyhedron.
In the present note we shall compute
the motivic nearby cycles $\cS_{P ({\mathbf f})}$
%%%%%%%%%%
%on the open
%set $U=X\setminus X_0(F)$
%%%%%%%%%%
of the composed function
$P(\mathbf{f})$
on  $X_0(\mathbf{f})$
as a sum over the set of compact
faces $\delta$ of the Newton polyedron of $P$.
For every such $\delta$,
let us denote by $P_\delta$ the corresponding quasi-homogeneous
polynomial. We associate to such
a quasi-homogeneous
polynomial a convolution operator
$\Psi_{P_\delta}$, which in the special
case
where $P_{\delta}$ is the polynomial $\Sigma = y_1 + y_2$
is nothing but the operator $\Psi_{\Sigma}$ considered in
\cite{ivc}.
For such a compact face $\delta$, one may also
define generalized
nearby cycles
$\cS_\mathbf{f}^{\sigma(\delta)}$, 
%%%fl
constructed
%%%%%%
as the limit, as $T \mapsto \infty$, of certain
truncated motivic zeta functions.

Our main result, Theorem \ref{MT}, follows from additivity from
the following statement, Theorem \ref{next}:
\begin{equation}\label{SfP*}
i^* \cS_{P(\mathbf{f}), U}=\sum_{\delta \in \Gamma^{\emptyset}}
 \Psi_{P_\delta}(\cS_\mathbf{f}^{\sigma(\delta)
}).
\end{equation}
Here $U$ denotes the complement of the locus 
where at least one function $f_j$ vanishes,
%%%%%%%%%mm
$\Gamma^{\emptyset}$
%%%%%%%%%%%%%
denotes
the set of compact faces of the Newton polyhedron of $P$
not contained in any coordinate hyperplane,
$\cS_{P(\mathbf{f}), U}$ refers to the extension of
$\cS_{P(\mathbf{f})}$ constructed in \cite{bittner2} and
\cite{ivc}, and $i^*$ denotes   restriction to $X_0
(\mathbf{f})$.

When $p = 2$ and $P = \Sigma$, one recovers
the motivic Thom-Sebastiani, cf. \cite{ts}, \cite{looi}, and \cite{barc},
in the way stated in \cite{ivc}.
When
$\mathbf{f}$ is the set of coordinate
functions on the affine space $\AA^p_k$,
%%fl
our result is equivalent to
 recovers a result obtained by
Guibert in \cite{guibert}.

This paper is a natural continuation
of \cite{ivc}, from which
part of the notation and several results are borrowed.

\section{Preliminaries}

\subsection{Grothendieck rings} \label{676} Throughout the
paper $k$ will be a field of characteristic zero. By a
variety over $k$, we mean a separated and reduced scheme of
finite type over $k$. If a linear algebraic group $G$ acts on a
variety $X$, we say the action is good if every $G$-orbit
is contained in an affine open subset of $X$. We denote by
$\Var^{G,\mathrm{eq}}$ the category of varieties with good
$G$-action, morphisms being $G$-equivariant morphisms. If
$S$ is a variety with good $G$-action, we denote by
$\Var_S^{G,\mathrm{eq}}$ the category of objects over $S$,
that is the category whose objects are morphisms
$Y\longrightarrow S$ in $\Var^{G,\mathrm{eq}}$, morphisms
in $\Var^{G,\mathrm{eq}}$ being defined in the standard
way.
Let $Y$ be a variety over $k$ and let $p : A \rightarrow Y$ be an affine bundle for the Zariski topology
(the fibers of $p$ are affine spaces and the transition morphisms between
trivializing charts are affine). In particular the fibers of $p$ have the structure of affine spaces.
Let $G$ be a linear  algebraic group. A good action of $G$ on $A$ is said to be affine if it is
a lifting of a good-action on $Y$ and its restriction to all fibers is affine.

One defines $K_0 (\Var_S^{G,\mathrm{eq}})$ as the free
abelian group on isomorphism classes of objects
$Y\longrightarrow S$ in $\Var_S^{G,\mathrm{eq}}$, modulo
the relations
\begin{equation} [Y \rightarrow S] = [Y' \rightarrow S] + [Y \setminus Y' \rightarrow S]
\end{equation}
for $Y'$ closed $G$-invariant in $Y$
and, for $f : Y \rightarrow S$ in
$\Var_{S}^{G,\mathrm{eq}}$,
 \begin{equation}
 [Y \times \AA^n_k \rightarrow
S, \sigma] = [Y \times \AA^n_k \rightarrow S, \sigma']
\end{equation} if $\sigma$ and $\sigma'$ are two liftings of the same
$G$-action on $Y$ to an affine action, the morphism $Y \times \AA^n_k
\rightarrow S$ being composition of $f$ with projection on
the first factor. Fiber product over $S$ induces a product
in the category $\Var_{S}^{G,\mathrm{eq}}$, which allows to
endow $K_0 (\Var_S^{G,\mathrm{eq}})$ with a natural ring
structure. Note that the unit $1_{S}$ for the product is
the class of the identity morphism $S\longrightarrow S$.

\subsection{} Let $s$ denote a positive integer and let $S$ be a
$k$-variety. From now on, we will consider only
$\GM^s$-actions on $S\times \GM^r$
which are trivial on the first factor.

Let us consider the category $\cC$ whose objects are
finite morphisms $\varphi : \GM^s \longrightarrow
\GM^{s'}$, a morphism between $\varphi : \GM^s
\longrightarrow \GM^{s'}$ and $\varphi ' : \GM^s
\longrightarrow \GM^{s''}$ being a finite morphism
$\vartheta
: \GM^{s'} \longrightarrow \GM^{s''}$ such that
$\vartheta\circ \varphi=\varphi '$.

We consider also the full subcategory $\cC '$ of $\cC$ the
objects of which are finite morphisms $\varphi : \GM^s
\longrightarrow \GM^{s}$. The subcategory $\cC'$ is final
in  $\cC$ in the language of \cite{ml}.

A morphism $\varphi : \GM^s \longrightarrow \GM^{s'}$
induces a natural functor
\begin{equation}
\Phi :\Var_{S\times \GM^r}^{\GM^{s'},\mathrm{eq}}
\longrightarrow \Var_{S\times \GM^r}^{\GM^{s},\mathrm{eq}}
\end{equation}
where an object $Y\longrightarrow {S\times \GM^r}$ with a
good $\GM^{s'}$-action is sent on the same underlying
object of $\Var_{S\times \GM^r}$ with the $\GM^{s}$-action
induced via $\varphi$.

The functor $\Phi$ induces a morphism
\begin{equation}\label{}
 K_0(\varphi)
:K_0(\Var_{S\times \GM^r}^{\GM^{s'},\mathrm{eq}})
\longrightarrow K_0(\Var_{S\times
\GM^r}^{\GM^{s},\mathrm{eq}}).
\end{equation}
We will denote by $K_0(\Var_{S\times
\GM^r}^{\varphi,\mathrm{eq}})$ the image of the morphism
$K_0(\varphi)$.

For every morphism $\vartheta$ between $\varphi$ and
$\varphi '$  in $\cC$, we get a morphism
\begin{equation}
K_0(\vartheta) :K_0(\Var_{S\times \GM^r}^{\varphi '
,\mathrm{eq}}) \longrightarrow K_0( \Var_{S\times
\GM^r}^{\varphi ,\mathrm{eq}})
\end{equation}
where a class of a good $\GM^s$-action induced by a
$\GM^{s''}$-action via $\varphi '$ on an object of
$\Var_{S\times \GM^r}$ is sent on the class of the same
$\GM^s$-action as induced by a $\GM^{s'}$-action via
$\varphi$. As a particular case, taking
$\varphi=\mathrm{Id}$, we get the natural inclusion of
$K_0(\Var_{S\times \GM^r}^{\varphi,\mathrm{eq}})$ into
$K_0(\Var_{S\times \GM^r}^{\GM^s,\mathrm{eq}})$.

We define the Grothendieck ring $K_0 (\Var_{S\times
\GM^r}^{\GM^s})$ as the colimit along $\cC$ (or along
$\cC'$, which amounts to the same) of the rings $K_0
(\Var_{S\times \GM^r}^{\varphi,\mathrm{eq}})$.

Note that we could also have defined
 the rings $K_0
(\Var_{S\times \GM^r}^{\varphi,\mathrm{eq}})$
and
$K_0 (\Var_{S\times
\GM^r}^{\GM^s})$ as suitable Grothendieck  rings
of the essential image
$\Var_{S\times
\GM^r}^{\varphi,\mathrm{eq}}$ of $\Phi$ and of
 the colimit
$\Var_{S\times \GM^r}^{\GM^s}$ along $\cC$
(or $\cC'$) of the categories $\Var_{S\times
\GM^r}^{\varphi,\mathrm{eq}}$, respectively.

There is a natural structure of $K_0 (\Var_k)$-module on
$K_0 (\Var_{S \times \GM^r}^{\GM^s})$. We denote by
$\LL_{S\times \GM^r} = \LL$ the element $\LL \cdot
1_{S\times \GM^r}$ in this module, and we set
\begin{equation}\cM_{S
\times \GM^r}^{\GM^s}: = K_0 (\Var_{S \times
\GM^r}^{\GM^s}) [\LL^{-1}].
\end{equation}
%Finally we denote by $\cN_{S
%\times \GM^r}^{\GM^s}$ the localization of the ring $\cM_{S
%\times \GM^r}^{\GM^s}$ with respect to the multiplicative
%family $\Bigl(\frac{1}{1 - \LL^{-i}}\Bigr)_{i > 0}$.

Note that when $s = r$ the above definitions $ K_0 (\Var_{S \times
\GM^r}^{\GM^s}) $
and $\cM_{S
\times \GM^r}^{\GM^s}$
coincides with that of \cite{ivc}
by \cite{ivc} 2.7.

A morphism $\vartheta :\GM^s\longrightarrow \GM^{s'}$
induces a morphism from $\cM_{S \times \GM^r}^{\GM^{s'}}$
to $\cM_{S \times \GM^r}^{\GM^s}$. For example, the
diagonal morphism $\GM \longrightarrow\GM^r$ yields a
canonical morphism
\begin{equation}\label{diag}
\Delta : \cM_{S \times \GM^r}^{\GM^{r}}\longrightarrow
\cM_{S \times \GM^r}^{\GM}
\end{equation}
Through this morphism, the class of a $\GM^{r}$-action
$\alpha$ on an object of $\Var_{S \times \GM^r}$ is sent on
the class of $\GM$-actions induced by $\alpha$ via a finite
group morphism from $\GM$ to $\GM^{r}$.

If $f : S \rightarrow S'$ is a morphism of varieties, composition with $f$ leads to a
push-forward morphism
$f_! : \cM_{S \times \GM^r}^{\GM^s} \rightarrow \cM_{S' \times \GM^r}^{\GM^s}$,
while fiber product
leads to a
pull-back morphism
$f^*: \cM_{S' \times \GM^r}^{\GM^s} \rightarrow \cM_{S \times \GM^r}^{\GM^s}$.

\subsection{}\label{sr}
Let $A$ be one of the rings
$\ZZ [\LL, \LL^{-1}]$,
$\ZZ [\LL, \LL^{-1}, (\frac{1}{1 - \LL^{-i}})_{i > 0}]$,
$\cM_{S \times \GM^r}^{\GM}$, etc.
%$\cN_{S \times \GM^r}^{\GM}$, etc.
We denote by
$A[[T]]_{\rm sr}$
the
$A$-submodule of
$A [[T]]$
generated by $1$ and by finite sums  of products of terms
$p_{e, i} (T) = \frac{\LL^e T^i}{1 - \LL^e T^i}$,
with $e$ in $\ZZ$ and $i$ in $\NN_{>0}$.
There is a unique
$A$-linear morphism
\begin{equation}
\lim_{T \mapsto \infty} : A[[T]]_{\rm sr} \longrightarrow A
\end{equation}
such that
\begin{equation}
\lim_{T \mapsto \infty} (\prod_{i \in I} p_{e_i, j_i} (T)) = (- 1)^{|I|},
\end{equation}
for every family $((e_i, j_i))_{i \in I}$ in
$\ZZ \times \NN_{>0}$,
with $I$ finite, maybe empty.

\subsection{}\label{ac}
We denote as usual by  $\cL_n (X)$
the space of arcs of order $n$,  also
known as the  $n$-th jet space on $X$.
It is a $k$-scheme whose $K$-points, for $K$ a field containing $k$,
is the set of morphisms
$\varphi : \Spec K [t]/ t^{n + 1} \rightarrow X$.
There are canonical morphisms $\cL_{n + 1} (X) \rightarrow
\cL (X)$ and the
arc space $\cL (X)$ is
defined as the projective limit of this system.
We denote by $\pi_n : \cL (X) \rightarrow \cL_n (X)$
the canonical morphism.
There is a canonical $\GM$-action on $\cL_n (X)$ and on $\cL (X)$
given by
$a \cdot  \varphi (t) = \varphi (at)$.

Let $X$ be a smooth variety over $k$ of pure dimension $d$ and
$g : X \rightarrow \AA^1_k$.
Set $X_0 (g)$ for the zero locus of $g$,
and define, for $n \geq 1$,
the variety
\begin{equation}\label{jkk}
\cX_n (g) := \Bigl \{\varphi \in \cL_n (X)  \Bigm  \vert \ord_t g
(\varphi) = n \Bigr\}.
\end{equation}
Note that $\cX_n (g)$ is invariant by the $\GM$-action on $\cL_n (X)$ and that
furthermore $g$ induces a morphism
$g_n : \cX_n (g) \rightarrow \GM$, assigning to a point
$\varphi$ in $\cL_n (X)$
the coefficient of $t^n$ in $g (\varphi)$, which we shall denote by
$\ac (g) (\varphi)$. We have
$g_n (a \cdot \varphi) = a^n g_n (\varphi)$,
hence with the terminology of \cite{ivc}
$g_n$ is diagonally monomial of weight $n$
with respect to the $\GM$-action
on $\cX_n (g)$. In particular, we may
consider the class $[\cX_n (g)]$ of
$\cX_n (g)$ in $\cM_{X_0 (g) \times \GM}^{\GM}$
and
the motivic zeta function
\begin{equation}
 Z_g(T): =\sum_{n\geq 1}[\cX_n(g)]\,\mathbf{L}^{-nd}\,T^n
\end{equation}
in $\cM_{X_0 (g) \times \GM}^{\GM}[[T]]$.

Denef and Loeser showed in \cite{motivic} and \cite{barc},
see also \cite{looi}
and \cite{ivc},
that $Z_g(T)$ is a rational series in
$\cM_{S \times \GM}^{\GM} [[T]]_{\rm sr}$
by giving a formula for
$Z_g(T)$ in terms of a resolution of $f$ we shall recall in
\ref{resolutions}.

\subsection{Resolutions}\label{resolutions}
Let us introduce some notation and terminology. Let $X$ be
a smooth variety of pure dimension $d$ and let $F$ a closed
subset of $X$ of codimension everywhere $\geq 1$. By a
log-resolution $h : Y \rightarrow X$ of $(X, F)$, we mean a
proper morphism $h : Y \rightarrow X$ with $Y$ smooth such
that the restriction of $h : Y \setminus h^{-1} (F)
\rightarrow X \setminus F$ is an isomorphism, and $h^{-1}
(F)$ is a divisor with normal crossings. We denote by
$E_i$, $i$ in $A$, the set of irreducible components of the
divisor $h^{-1}  (F)$. For $I \subset A$, we set
\begin{equation}E_I :=\bigcap_{i \in I} E_i
\end{equation}and
\begin{equation}
E_I^{\circ} := E_I \setminus \bigcup_{j \notin I} E_j.
\end{equation}
We denote by $\nu_{E_i}$  the normal bundle of $E_i$ in
$Y$ and
by $\nu_{E_I}$ the fiber product
of the restrictions to $E_I$ of the bundles
$\nu_{E_i}$, $i$ in $I$.
We will denote by $U_{E_i}$ the complement of the zero
section in $\nu_{E_i}$ and by $U_I$ the fiber product
of the restrictions of the spaces
$U_{E_i}$, $i$ in $I$, to $E_I^{\circ}$.

\bigskip
If $\cI$ is an ideal sheaf  defining a closed subscheme of
$X$
and $h^{*} (\cI)$ is locally principal, we define $N_i
(\cI)$, the multiplicity of $\cI$ along $E_i$, by the
equality of divisors
\begin{equation}
h^{-1} (F) = \sum_{i \in A} N_i (\cI) E_i.
\end{equation}
If $\cI$ is principal generated by a function $g$ we write
$N_i (g)$ for $N_i (\cI)$. Similarly,  we define integers
$\nu_i$ by the equality of divisors
\begin{equation}
K_Y = h^*{K_X} + \sum_{i \in A} (\nu_i - 1) E_i.
\end{equation}

\subsection{}\label{g_I}
Assume again $g$ is a function on a smooth variety $X$ of
pure dimension $d$. Let $F$ a reduced divisor containing
$X_0 (g)$ and let $h : Y \rightarrow X$ be a log-resolution
of $(X, F)$. Let us explain how $g$ induces a morphism $g_I
: U_I \rightarrow \GM$. Note that the function $g \circ h$
induces a function
\begin{equation}\label{gtt}
\bigotimes_{i \in I} \nu_{E_i}^{\otimes N_i (g)}{}_{|E_I}
\longrightarrow \AA^1_k,
\end{equation}
vanishing only on the zero section. We define $g_I
: \nu_{E_I} \rightarrow \AA^1_k$ as the composition of this
last function with the natural morphism $\nu_{E_I}
\rightarrow \otimes_{i \in I} \nu_{E_i}^{\otimes N_i
(g)}{}_{|E_I}$, sending $(u_i)$ to $\otimes u_i^{\otimes
N_i (g)}$. We still denote by $g_I$ the induced morphism
from $U_I$ (resp. $U_{E_I}$) to $\GM$.

We view $U_I$ as a variety over $X_0 (g) \times \GM$ via
the morphism $(h\circ\pi_I,g_I)$. The group $\GM$ has a
natural action on each $U_{E_i}$, so the diagonal action
induces a $\GM$-action on $U_I$. Furthermore, the morphism
$g_I$ is monomial, in the terminology of \cite{ivc}, hence $ U_I \rightarrow X_0 (g) \times
\GM$ has a class in $\cM_{X_0 (g) \times \GM}^{\GM}$ which
we will denote by $[U_I]$.

\subsection{}We now assume that $F = X_0 (g)$, that
is $h : Y \rightarrow X$ is a log-resolution of $(X,
X_0(g))$. In this case, $h$ induces a bijection between
$\cL (Y) \setminus \cL (\vert h^{-1} (X_0 (g))\vert)$ and
$\cL (X) \setminus \cL (X_0 (g))$.

One deduces from Lemma 3.4 in \cite{arcs}, in a way completely similar to
\cite{motivic} and \cite{barc}, the equality
\begin{equation}\label{dl1}
Z_g (T)  = \sum_{\emptyset \not= I \subset A} [U_I]
\prod_{i \in I} \frac{1}{T^{- N_i (g)} \LL^{\nu_i} - 1}
\end{equation}
in $\cM_{X_0 (g) \times \GM}^{\GM}[[T]]$.

In particular, the function  $Z_g (T)$ is rational and
belongs to $\cM_{X_0 (g) \times \GM}^{\GM} [[T]]_{\rm sr}$,
with the notation of \ref{sr}, hence we can consider $
\lim_{T \mapsto \infty} Z_g (T)$ in $\cM_{X_0 (g) \times
\GM}^{\GM}$
and
set
\begin{equation}
\cS_g := - \lim_{T \mapsto \infty} Z_g (T),
\end{equation}
which by (\ref{dl1})  may be expressed on a resolution $h$
as
\begin{equation}\label{dl2}
\cS_g  = - \sum_{\emptyset \not=  I \subset A}
(-1)^{|I|}[U_I],
\end{equation}
in $\cM_{X_0 (g) \times \GM}^{\GM}$.
The element $\cS_g$ is called
the motivic Milnor fiber or the motivic nearby fiber
of $f$. It
was first considered by Denef and Loeser,
cf.  \cite{motivic}, \cite{barc} and \cite{lef}.
For recent results concerning $\cS_g$,
we refer the reader to \cite{guibert},
\cite{bittner2} and \cite{ivc}.

\subsection{}
Consider a family $\mathbf{f}$ of $p$  functions $f_j : X
\rightarrow \AA^1_k$, $1\leq j\leq p$. We denote by
$X_0(\mathbf{f})$ the set of common zeroes of the functions
$f_j$, $1\leq j\leq p$ and by $F$ the product function
$f_1\ldots f_p$.

Let us fix a rational polyhedral
convex cone $C$ in $\RR_{>0}^p$ and
an integral linear form $\ell$ on $\ZZ^p$ which is positive
on $\bar C\setminus \{0\}$,
%%%%%%%%%%%%%%%fl
where $\bar C$ denotes the closure of $C$ in $\RR^p$.
%%%%%%%%%%%%

We will consider the modified zeta function
$Z^{C,\ell}_{\mathbf{f}}$ defined as follows: for a vector
$\mathbf{n}$ in $\NN^p_{>0}$, we denote by $s(\mathbf{n})$
the sum of its components and we consider, similarly as in
(\ref{jkk}), the variety
\begin{equation}\cX_\mathbf{n} (\mathbf{f})
:= \Bigl \{\varphi \in \cL_{s(\mathbf{n})} (X)  \Bigm  \vert
\ord f_j(\varphi) = n_j
%%%%%%%%%%%%%%%mm
,\,1\leq j\leq p
%%%%%%%%%%%%%%%
\Bigr\}.
\end{equation}
Note that
$\cX_n (\mathbf{f})$ is
%%%%%%%%%%%%%%%mm fl
stable under
%%%%%%%%%%%%%%%
the $\GM$-action on $\cL_n (X)$ and that $\mathbf{f}$
induces a morphism
\begin{equation}
\mathbf{f_n} :
\cX_n (\mathbf{f})
\longrightarrow
\GM^p,
\end{equation}
whose components are
%%%%%%%%%%%%%mm
$\ac (f_j),\,1\leq j\leq p$
%%%%%%%%%%%%
defined similarly as in \ref{ac}. Since $\mathbf{f_n} (a 
%%%fl
\cdot
%%%%
\varphi) = a^{\mathbf{n}} \mathbf{f_n} (\varphi)$, we may
consider the class $[\cX_\mathbf{n} (\mathbf{f})]$ of
$\cX_n (\mathbf{f}) \rightarrow X_0 (\mathbf{f}) \times
\GM^p$  in $\cM_{X_0 (\mathbf{f}) \times \GM^p}^{\GM}$. We
set
\begin{equation}Z^{C,\ell}_{\mathbf{f}} (\mathbf{T})
:= \sum_{\mathbf{n}\in C} [\cX_\mathbf{n}
(\mathbf{f})]\, \mathbf{L}^{-s(\mathbf{n})
d}\,T^{\ell(\mathbf{n})}
\end{equation} in $\cM_{X_0 (\mathbf{f})
\times \GM^p}^{\GM} [[T]]$.

\subsection{} Let $h : Y \rightarrow X$ be a log-resolution of the
set  $X_0(F)$. We keep the notations of
\ref{resolutions}. In particular we denote by $A$ the set of
irreducible components of $h^{-1}(X_0(F))$. For $i$ in $A$ we
will denote by $N_i$ the integral vector of the orders
$N_i(f_j)$ of the functions
$f_j$, $1\leq j\leq p$, along the divisor
$E_i$, and by $N_I$ the linear map
\begin{equation}
N_I : \begin{cases}\RR_{>0}^I &\longrightarrow \RR_{>0}^p\\
\mathbf{k}&\longmapsto
\sum_{i\in I} k_i N_i.
\end{cases}
\end{equation}
Similarly the set of integers $\nu_i$ defines a linear integral form
$\nu_I : \mathbf{k}\longmapsto \sum_{i\in I} k_i \nu_i$ on
$\RR_{>0}^I$.

Using Lemma 3.4 in \cite{arcs} similarly as for the proof of (\ref{dl1}),
see for example \cite{barc}, \cite{looi},
one gets the following formula for
the zeta function $Z^{C,\ell}_{\mathbf{f}} (T)$
in terms of the resolution:
\begin{equation}
Z^{C,\ell}_{\mathbf{f}} (T)= \sum_{\emptyset\neq I\subset
A}[U_I]\sum_{\{\mathbf{k}\in
\NN_{>0}^p \vert N_I (\mathbf{k}) \in C\}}\prod_{i\in I} (T^{\ell(N_i)}
\LL^{-\nu_i})^{k_i}.
\end{equation}
Here $[U_I]$
stands for the class in $\cM_{X_0 (\mathbf{f}) \times
\GM^p}^{\GM}$ of the morphism $(h,\mathbf{f}_I) :
U_I\longrightarrow X_0 (\mathbf{f}) \times\GM^p$.

It  follows that $Z^{C,\ell}_{\mathbf{f}} (T)$
belongs to $\cM_{X_0 (\mathbf{f}) \times
\GM^p}^{\GM}[[T]]_{\rm sr}$, hence we may set
\begin{equation}\label{hhh}
\cS_{\mathbf{f}}^{C,\ell} := \lim_{T \mapsto \infty}
Z^{C,\ell}_{\mathbf{f}} (T)
\end{equation}
in $\cM_{X_0 (\mathbf{f}) \times \GM^p}^{\GM}$. By
section 2.9 of \cite{ivc}
we have:
\begin{equation}\label{SC}
\cS_{\mathbf{f}}^{C,\ell} = \sum_{I}\chi(N_I^{-1}(C))[U_I],
\end{equation}
where $\chi$ denotes Euler characteristic with compact
supports. Note that this is independent of $\ell$, so we
%%fl
may 
write $\cS_{\mathbf{f}}^{C}$ instead of
$\cS_{\mathbf{f}}^{C,\ell}$.

\section{Composition with a non-degenerate polynomial}

\subsection{The generalized convolution  $\Psi_P$}

Let $P$ be a quasi-homogeneous polynomial function on
$\GM^p$, that is $P$ is homogeneous for a $\GM$-action
$\alpha$ on $\GM^p$ monomial of weight
$\mathrm{w}=(w_1,\ldots,w_p)$.

Let $X$ be a  smooth variety. We will denote by $\pr_1$ the
projection of $X\times \GM^p\times \GM$ on  $X\times \GM$
(forgetting the $\GM^p$ factor) and by $i$ the inclusion of
the complement of $X\times P^{-1}(0)$ into $X\times \GM^p$.

For a variety $A$ of dimension $e$ in
$\mathrm{Var}_{X\times \GM^p}$, the function $P$ induces by
composition with the second projection a function on $A$ we
still denote by $P$
\begin{equation}
 P :   A \longrightarrow \AA^1_k.
\end{equation}

We now define the (augmented) zeta function $Z^0_P(T)$ as
\begin{equation}\label{zeta0}
Z^0_P(T)=\sum_{n\geq
0}[\cX_n(P)]\LL^{-ne}T^n=[\cX_0(P)]+Z_P(T),
\end{equation}
where $\cX_n(P)$ is
\begin{equation}
\cX_n (P) := \Bigl \{\varphi \in \cL_n (A)  \Bigm  \vert
\ord_t P (\varphi) = n \Bigr\},
\end{equation}
for $n\geq 0$.
It  belongs to
$\cM_{X\times \GM^p \times \GM}^{\GM} [[T]]_{\rm sr}$.
We define
$\Psi^0_P(A)$ as the limit, as $T \mapsto \infty$,
of the opposite
$-Z^0_P(T)$. Thus, with the notations of
\cite{ivc},  it  is nothing but
\begin{equation}
-\lim_{T\rightarrow\infty}Z^0_P(T)=-[A\setminus
P^{-1}(0)]+\cS_P([A]).
\end{equation}
It is an object in $\cM_{X\times \GM^p \times \GM}^{\GM}$,
the $\GM$-action and the morphism to  $\GM$
being the usual
ones. On $A\setminus P^{-1}(0)$ the $\GM$-action is
trivial and the morphism to $\GM$ is the restriction of $P$
to $A\setminus P^{-1}(0)$. Taking the direct image by the
projection $\pr_1$ we get the following
%%%%%%%%%%%%%mm
%well defined
%%%%%%%%%%%%
object in
$\cM_{X\times \GM}^{\GM}$:
\begin{equation}\label{psiA}
\Psi^0_P(A):={\pr_1}_!\,(-[A\setminus P^{-1}(0)]+\cS_P(A)).
\end{equation}
One may then extend uniquely this construction to a
$\cM_k$-linear group morphism
\begin{equation}
\Psi_P^0 : \cM_{X\times \GM^p} \longrightarrow \cM_{X\times
\GM}^{\GM}.
\end{equation}

If $A$ is endowed with a $\GM$-action $\alpha$ for which
the morphism to $\GM^p$ is monomial of weight $\mathrm{w}$,
$A\setminus P^{-1}(0)$ is endowed with an additional action
which is homogeneous with respect to the composed morphism
to $\GM$. Hence we may attach to $A\setminus P^{-1}(0)$ a
class $[A\setminus P^{-1}(0)]$ in $\cM_{X\times \GM^p
\times \GM}^{\GM^2}$. In
%%%%%%%%%%%%%%mm
\cite{ivc}, \S \kern .15em 3.10,
%%%%%%%%%%%%
we attached
to such an  $A$ with the  action $\alpha$ an element
$\cS_P(A)$
in $\cM_{X\times \GM^p \times \GM}^{\GM^2}$.
Hence we can consider
${\pr_1}_!\,(-[A\setminus P^{-1}(0)]+\cS_P(A))$
as an element of $\cM_{X \times
\GM}^{\GM^2}$. Composing with the canonical morphism
$\cM_{X \times
\GM}^{\GM^2} \rightarrow
\cM_{X \times
\GM}^{\GM}$
%%%%%%%%%%%%mm
induced by
%%%%%%%%%%
the diagonal action, we get an element of $\cM_{X \times
\GM}^{\GM}$ we shall denote by $\Psi_P(A)$. This
construction extends uniquely to a $\cM_k$-linear group
morphism
\begin{equation}
\Psi_P: \cM_{X\times \GM^p}^{\GM}\longrightarrow
\cM_{X\times \GM}^{\GM}.
\end{equation}

%\begin{prop}Let $X$ be a smooth variety and $P$ a $p$-variable quasi
%homogeneous polynomial. Then there exists a group morphism
%defined by
%\begin{equation}\label{PsiP}
%\Psi_{P} : \begin{cases}\cM_{X \times \GM^p}^{\GM} &\longrightarrow
%\cM_{X \times \GM}^{\GM}\\
 %A&\longmapsto {\pr_1}_!\,(-i^*(A)+\cS_P(A)).
%\end{cases}
%\end{equation}
%\end{prop}

\begin{remark}
When $P$ is the sum of coordinates $\Sigma$ on $\GM^2$,
then $\Psi_\Sigma$ is nothing but the convolution product
from \cite{ivc}. More precisely, the convolution product
$\Psi_\Sigma$ defined in \cite{ivc} is equal to the
composition of the morphism $\Psi_\Sigma$ defined in this
paper with the morphism $\Delta$ defined in (\ref{diag}).
\end{remark}

\subsection{Composed maps}
 For $1\leq j\leq p$, let
$f_j:X_j \longrightarrow \AA_k^1$ be a function on a smooth $k$-variety $X_j$.
By composition with the projection, $f_j$ becomes a
function on the product $X=\prod_j X_j$.
We write $d$ for the dimension of $X$. Define
$\mathbf{f}$ as the family of the $f_j$ on $X$, $1\leq j\leq
p$. The product of the log-resolutions of the
$X_{j,0}(f_j)$ is a log-resolution $h:Y\longrightarrow X$
of $X_0(F)$ (recall that $F=f_1\ldots f_p$).

Let $P= \sum_{\alpha \in \NN^p} a_{\alpha} y^{\alpha}$
be a polynomial in
$k[y_1,\ldots,y_p]$.
We denote by ${\rm supp} (P)$ the set of exponents $\alpha$ in $\NN^p$
with $a_{\alpha} \not= 0$.
The Newton polyhedron
$\Gamma$ of $P$ is the convex hull of
${\rm supp} (P) + \RR_+^p$.
For a compact face $\delta$ of $\Gamma$
we  denote by $P_{\delta}$ the sum of the
monomials of $P$ supported in $\delta$:
\begin{equation}
P_{\delta} = \sum_{\alpha \in \delta} a_{\alpha} y^{\alpha}.
\end{equation}
We say
$P$ is
non-degenerate with respect to its Newton polyhedron $\Gamma$,
if,  for every compact face $\delta$ of
$\Gamma$, the function  $P_{\delta}$ is smooth on $\GM^p$.

To the Newton
%%%%%%%%%%%%mm
polyhedron
%%%%%%%%%%%
$\Gamma$ one may associate a fan
of rational
%%%%%%%%%%%%mm
polyhedral
%%%%%%%%%%
cones subdividing  $\RR_+^p$ as
follows. We consider the function $\ell_{\Gamma}$ assigning
to a vector $a$ in $\RR_+^p$ the value $\inf_{b \in \Gamma}
\langle a ,  b \rangle$, with $\langle ,  \rangle$ the
standard inner product. For any $a$ in $\RR_+^p$ we may
consider the compact face
\begin{equation}
\delta_a = \{b \in \Gamma_{c} \vert
\langle a ,  b \rangle = \ell_{\Gamma} (b) \},
\end{equation}
 with $\Gamma_{c}$ the union of all compact faces of $\Gamma$.

For a compact face $\delta$ of the Newton polyhedron $\Gamma$, we
denote by $\sigma(\delta)$ its dual cone
$\{a \in \RR_+^p \vert \delta_a = \delta \}$.
The cones
$\sigma(\delta)$, for
%%%%%%%%%mm
$\delta$
%%%%%%%%%
running over the compact faces of $\Gamma$,
form a fan
partitioning $\RR_+^p$ by rational polyhedral cones.
The function $\ell_{\Gamma}$ is linear on each cone
$\sigma(\delta)$.

We write
%%%%%%%%%%%%mm
$\Gamma_c $
%%%%%%%%%%%%
for the set of compact faces of $\Gamma$. For $J$ a subset
of $\{1, \dots, p\}$, we denote by
%%%%%%%%%%%%mm
$\Gamma^J $
%%%%%%%%%%%
the set of compact faces of $\Gamma$
contained in the coordinate hyperplanes $x_i = 0$ for $i$
in $J$, and in no other coordinate hyperplane, so that
%%%%%%%%%%%%mm
$\Gamma_c $
%%%%%%%%%%%%
is the disjoint union of the subsets
%%%%%%%%%%%%mm
$\Gamma^J $. Note that $\ell_\Gamma$ is positive on
$\overline{\sigma(\delta)}\setminus \{0\}$ if and only if
$\delta$ is in $\Gamma^\emptyset $.
%%%%%%%%%
We denote by $X_J$ the closed subset of $S$ defined by
the vanishing of the functions $f_i$, $i \in J$, and by
$\mathbf{f}_J : X_J \rightarrow \AA^{\{1, \dots, p\}
\setminus J}$ the morphism induced by the functions $f_j$,
$j \notin J$.

For every variety $Z$ containing $X_0 (\mathbf{f})$
we denote by $i^*$ the restriction morphisms
\begin{equation}
\cM_{Z \times \GM}^{\GM} \longrightarrow
\cM_{X_0 (\mathbf{f}) \times \GM}^{\GM}
\end{equation}
and
\begin{equation}
\cM_{Z \times \GM}^{\GM} [[T]]\longrightarrow
\cM_{X_0 (\mathbf{f}) \times \GM}^{\GM} [[T]].
\end{equation}

\begin{theorem}\label{MT}With the previous notations
and hypotheses, we have  the following formula  for
$i^* \cS_{P(\mathbf{f})}$ in $\cM_{X_0 (\mathbf{f}) \times \GM}^{\GM}$:
\begin{equation}\label{SfP}
i^*\cS_{P(\mathbf{f})}=\sum_{J \subset \{1, \dots, p\} }
\sum_{\delta \in
%%%%%%%%%%%mm
\Gamma^J
%%%%%%%%%%%
}
\Psi_{P_\delta}(\cS_{\mathbf{f}_J}^{\sigma(\delta),\ell_{\Gamma}}).
\end{equation}
\end{theorem}

\begin{proof}Following \cite{ivc}, for $\gamma$ in $\NN_{>0}$,
we consider the constructible set
\begin{equation}\cX_n^{\gamma n} :=
\Bigl\{\varphi \in
\cL_{\gamma n} (X)  \Bigm  \vert \ord_t
P(\mathbf{f})(\varphi) = n, \ord_t F(\varphi) \leq \gamma n
\Bigr\}
\end{equation}
together with the morphism
$\ac (P(\mathbf{f})) : \cX_n^{\gamma n} \rightarrow \GM$,
giving rise to a class
$[ \cX_n^{\gamma n}]$ in
$\cM^{\GM}_{X_0 (F) \times  \GM}$.
By Proposition 3.8 in
\cite{ivc}, for $\gamma \gg 0$,
the corresponding zeta function
\begin{equation}
Z^\gamma_{P(\mathbf{f}),X \setminus X_0 (F)}  (T)
:= \sum_{n>0}
 [ \cX_n^{\gamma n}]
\LL^{-\gamma nd} T^n
\end{equation}
lies in
$\cM^{\GM}_{X_0 (F) \times  \GM} [[T]]_{\rm sr}$
and its limit as $T \mapsto \infty$ is independent of $\gamma$, so we may set
\begin{equation}
\cS_{P(\mathbf{f}), X \setminus X_0 (F)} := - \lim_{T \mapsto \infty}
Z^\gamma_{P(\mathbf{f}),X \setminus X_0 (F)}  (T).
\end{equation}
Furthermore, by additivity of $\cS_{P(\mathbf{f})}$, cf. Theorem
%%%%%%%%%%mm
3.11 of \cite{ivc},
%%%%%%%%%%
we have
\begin{equation}
\cS_{P(\mathbf{f})} =\sum_{J \subset \{1, \dots, p\} }
\cS_{P(\mathbf{f})_{\vert X_J}, X_J^{\circ}},
\end{equation}
with $X_J^{\circ}$ the largest open in $X_J$, where no
$f_j$, $j \notin J$, vanishes.
Theorem \ref{MT} follows now directly from Theorem \ref{next}.
\end{proof}

\begin{theorem}\label{next}With the previous notation the following holds:
\begin{equation}\label{Spsi}
i^* \cS_{P(\mathbf{f}), X\setminus X_0(F)}=\sum_{\delta\in
%%%%%%%%%%%mm
\Gamma^J
%%%%%%%%%%%
}\Psi_{P_\delta}(\cS_\mathbf{f}^{\sigma(\delta)
%%%%%%%%%%%mm
%,\ell_{\Gamma}
%%%%%%%%%
}).
\end{equation}
\end{theorem}

\begin{proof}
Let us fix a log-resolution $h:Y\longrightarrow X$
of $X_0(F)$.  We shall
keep the notations of (\ref{resolutions}).

Fix a subset of $I$ of $A$ and
${\bf k} = (k_i)_{i \in I}$ in $\NN_{> 0}^{I}$.
For $\varphi$ in $\cL_{\gamma n} (Y)$ with $\varphi (0)$ in $E_i$,
we set $\ord_{E_i} \varphi := \ord_t z_i (\varphi)$, for
$z_i$ any local equation of $E_i$ at $\varphi (0)$.
We denote by
$\cX_{n,\mathbf{k}}$ the set of arcs
$\varphi$ in $\cL_{\gamma n} (Y)$ such that $\varphi(0)$ is in
$E_I^{\circ}$ and $\ord_{E_i} \varphi = k_i$ for $i \in I$.
We also consider the subset $\cY_{n, {\bf k}}$ of
$\cX_{n,\mathbf{k}}$ consisting of arcs $\varphi$ such that
$\ord_t ((P (\mathbf{f}))\circ h) (\varphi) = n$. The variety $\cY_{n,
{\bf k}}$ is stable by the usual $\GM$-action on
$\cL_{\gamma n}(Y)$ and the morphism $\ac(P (\mathbf{f}))\circ h)$ defines a class
$[\cY_{n, {\bf k}}]$ in $\cM_{X_0 (F) \times
\GM}^{\GM}$. Note that
$\cY_{n, {\bf k}} = \emptyset$
if $n <\ell_{\Gamma} (N_I (\mathbf{k}))$.

By a now standard calculation, using Lemma 3.6 in \cite{arcs},
$Z^\gamma_{P(\mathbf{f}),X \setminus X_0 (F)}$ may be
expressed on the log-resolution $Y$,
as
\begin{equation}\label{hotrats}
Z^{\gamma}_{P(\mathbf{f}),X \setminus
X_0(F)}=\sum_{\emptyset \not= I \subset A}
%%%%%%%%%%%%%%mm
\sum_{\sur{ \sur{N_I(\mathbf{k})\in \sigma
(\delta)}{n=\ell_{\Gamma} (N_I (\mathbf{k}))}} {\langle N_I
(\mathbf{k}), \bf{1} {\rangle} \leq \gamma n}}
%%%%%%%%%%%%
[\cY_{n, \mathbf{k}}] \, \LL^{\ \sum_{i \in I} (\nu_i - 1)
k_i} \LL^{- \gamma n d} T^n.
\end{equation}

We denote by $B$ the set of all subsets $I$ of $A$ such
that $h (E_I^{\circ})$ is contained in $X_0 (\mathbf{f})$.
We fix $I$ in $B$ and
${\bf k} = (k_i)_{i \in I}$ in $\NN_{> 0}^{I}$.
Note that there is a unique compact face $\delta$ of $\Gamma$
such that
$N_I ({\bf k}) $ lies in $\sigma (\delta)$.

To go further on,
we shall use the following variant
of the classical deformation to the normal cone already
considered in \cite{ivc}.
We consider the affine line $\AA^1_k = \Spec k
[u]$ and the subsheaf
\begin{equation}
\cA_{\bf k} := \sum_{{\bf n} \in \NN^I} \cO_{Y \times
\AA^1_k} \left(- \sum_{i \in I} n_i (E_i \times \AA^1_k) \right ) u^{-
\sum_{i \in I} k_i n_i}
\end{equation}
of $\cO_{Y \times \AA^1_k} [u^{-1}]$. It is a sheaf of
rings and we set
\begin{equation}
CY_{\bf k} := \Spec \cA_{\bf k}.
\end{equation}
The natural inclusion $\cO_{Y \times \AA^1_k} \rightarrow
\cA_{\bf k}$ induces a morphism $\pi : CY_{\bf k}
\rightarrow Y \times \AA^1_k$, hence a morphism $p :
CY_{\bf k} \rightarrow \AA^1_k$.
Via the same inclusion,
the functions $P ({\bf f}\circ h)$ and $F \circ h$
are, in $\cA_{\bf k}$, divisible by $ u^{\ell_{\Gamma} (N_I ({\bf k}))}$
and by $u^{\langle N_I ({\bf k}) , \mathbf{1} \rangle}$, where
$\mathbf{1}$ denotes the vector with all coordinates equal
to $1$, and we
denote the corresponding
quotients by $\widetilde P (\mathbf{f})_\mathbf{k}$
and $\widetilde
F_\mathbf{k}$, respectively.

We denote by $\widetilde E_i$ the strict transform of the divisor
$E_i \times \AA^1_k$ by $\pi$,
by
$D$ the divisor globally defined on $CY_\mathbf{k}$
by $u=0$, and by $CE_i$ the divisors $\widetilde E_i-k_iD$,
$i$ in $I$.
We denote by $CY_\mathbf{k}^\circ$ the
complement in $CY_\mathbf{k}$ of the union of the divisors $CE_i$
$i$ in $I$,
and we denote by $Y^\circ$ the complement in $Y$ of the
union of the $E_i$, $i$ in $I$.

As proved in
%%%%%%%%mm
Lemma 5.12 of \cite{ivc},
%%%%%%%%%%%%%
the scheme $CY_{\bf k}$ is
smooth, the morphism $\pi$ induces an isomorphism above
$\AA^1_k \setminus \{0\}$, the morphism $p$ is a smooth
morphism and its fiber $p^{-1} (0)$ may be naturally
identified with the bundle $\nu_{E_I}$. Furthermore, when
restricted to $CY_{\bf k}^\circ$, the fiber of $p$ above
$0$ is naturally identified with $U_I$ and $\pi$ induces an
isomorphism between
 $CY_{\bf k}^\circ \setminus p^{-1}(0)$ and $Y^\circ\times \AA^1_k
\setminus \{0\}$.
The restrictions of $\widetilde P (\mathbf{f})_\mathbf{k}$
and
$\widetilde
F_\mathbf{k}$
to the
fiber $U_I \subset p^{-1} (0)$ are respectively equal to
$P_{\delta} ({\mathbf f}_I)$ and
$F_I$.

The ring $\cA_\mathbf{k}$ being a graded subring of the ring
$\cO_Y [u, u^{-1}]$, we may consider
the $\GM$-action $\sigma$ on $CY_{\bf k}$,
leaving sections of $\cO_Y$ invariant and acting by
$\lambda \mapsto \lambda^{-1} u$ on $u$.
It restricts
on $U_I$, to
the diagonal action induced by the canonical $\GM^I$-action
on $U_I$ via the finite morphism $\lambda \longmapsto
\lambda^\mathbf{k}$.
We have now two different $\GM$-actions
on $\cL_n(CY_\mathbf{k}^\circ)$ : the one induced by the standard $\GM$-action
on arc spaces and the one induced by $\sigma$.
We denote by $\widetilde \sigma$ the action given by the composition of these
two (commuting) actions.

Let us denote by $\widetilde \cL_{\gamma n}(CY_\mathbf{k}^\circ)$ the
set of arcs $\varphi$ in $\cL_{\gamma n}(CY_\mathbf{k}^\circ)$ such
that $p(\varphi(t))=t$ (in particular $\varphi(0)$ is in $U_I$).
For such an arc $\varphi$, composition with $\pi$ send
$\varphi$ to an arc in $\cL_{\gamma n}(Y\times\AA^1_k)$ which is the
graph of an arc in $\cL_{\gamma n}(Y)$ not contained in the union of
the divisors $E_i$, $i$ in $I$.
Note that $\widetilde \cL_n(CY_\mathbf{k}^\circ)$ is stable by $\widetilde \sigma$.

\begin{lem}\label{rio}Let $I$ be in $B$ and $\mathbf{k}$ in $\NN_{> 0}^I$.
Assume $n \geq k_i$ for $i$ in $I$.
The morphism
$\widetilde \pi
: \widetilde \cL_n(CY_\mathbf{k}^\circ) \rightarrow \cX_{n,\mathbf{k}}$
induced by the projection $CY_{\bf k}^\circ \rightarrow Y$
is an affine bundle with fiber $\AA^{\sum_I k_i}_k$.
Furthermore if  $\widetilde \cL_n(CY_\mathbf{k}^\circ)$ is
endowed with the $\GM$-action induced by $\widetilde
\sigma$ and $\cX_{n,\mathbf{k}}$ with  the standard
$\GM$-action, $\widetilde \pi$ is $\GM$-equivariant and the
action of $\GM$ on the affine bundle is affine.
Furthermore, if $n \geq \ell_{\Gamma}(N_I (\mathbf{k}))$,
then for every $\varphi$ in $ \widetilde \cL_{\gamma
n}(CY_\mathbf{k}^\circ)$
\begin{equation}
\ac (P ({\bf{f}}) \circ h) (\widetilde \pi (\varphi))=\ac
(\widetilde P({\bf f})_{\bf k} (\varphi)).
\end{equation}
When $P_{\delta} ({\bf f}_I)
(\varphi(0))\neq 0$, hence $(\ord_t(P ({\bf f})) \circ h)(\widetilde \pi (\varphi)) 
= \ell_{\Gamma} (N_I ({\bf k}))$,
 we have
\begin{equation}\ac (P ({\bf f})) \circ h) (\widetilde \pi
(\varphi))= P_{\delta} ({\bf f}_I) (\varphi(0)).
\end{equation}
\end{lem}

\begin{proof}The first part  statement is contained in
%%%%%%%%mm
Lemma 5.13 of \cite{ivc}
%%%%%%%%%%%%%
and the rest follows from its proof.
\end{proof}

We then define $\widetilde \cY_{n, {\bf k}}$ as the inverse
image of $\cY_{n, {\bf k}}$ by the fibration $\widetilde\pi$.
It is the subset of arcs $\varphi$ in
$\cL_{\gamma n}(CY_{\mathbf{k}}^\circ)$ such that
$\ord_t \widetilde
P(\mathbf{f})_{\mathbf{k}}(\varphi)=n- \ell_{\Gamma} (N_I (\mathbf{f}))$.
We denote
by $[\widetilde \cY_{n, {\bf k}}]$ the class of $\widetilde \cY_{n,
{\bf k}}$ in $\cM_{X_0 (F) \times
\GM}^{\GM}$, the morphism $\widetilde \cY_{n, {\bf k}}
\rightarrow \GM$ being $\ac (\widetilde P(\mathbf{f})_{\bf k} )$ and the
$\GM$-action being induced
by $\widetilde \sigma$.
We denote by $[U_I \setminus
(P_{\delta} (\mathbf{f}_I)^{-1} (0))]$ the class of
$U_I \setminus (P_{\delta} (\mathbf{f}_I)^{-1} (0))$ in
$\cM_{X_0 (F) \times \GM}^{\GM}$,
the $\GM$-action being the natural diagonal action of
weight ${\bf k}$ on $U_I \setminus
(P_{\delta} (\mathbf{f}_I)^{-1} (0))$ and the
morphism to $\GM$ being the restriction of $P_{\delta} (\mathbf{f}_I)$.
We also
consider the class $[\GM \times F_I^{-1} (0)]$ of $\GM
\times P_{\delta} (\mathbf{f}_I)^{-1} (0)$ in $\cM_{X_0 (F)
\times \GM}^{\GM}$, the $\GM$-action on the second factor
being the diagonal one and the morphism to $\GM$ being the
first projection.

\begin{lem}\label{cal}Let $I$ be in %%%j'insiste !!
 $B$ and $\mathbf{k}$ in $\NN_{> 0}^I$.
The following equalities
hold  in $\cM_{X_0 (F) \times \GM}^{\GM}$
\begin{enumerate}
\item[(1)]
$[\widetilde \cY_{n, {\bf k}}] =\LL^{\gamma nd}
[U_I \setminus (P_{\delta} (\mathbf{f}_I)^{-1} (0))]$,
if $n = \ell_{\Gamma} (N_I (\mathbf{k}))$,
\item[(2)]$[\widetilde \cY_{n, {\bf k}}]= \LL^{\gamma nd -m} [\GM \times P_{\delta} (\mathbf{f}_I)^{-1} (0)]$,
if $n - \ell_{\Gamma} (N_I (\mathbf{k}))= m> 0$.
\end{enumerate}
\end{lem}

\begin{proof}
As we assume $P$ is non-degenerate with respect to its
Newton polyhedron, $P_\delta$ is smooth on $\GM^p$ and the
composed map $P_\delta(\mathbf{f}_I)$ is smooth on $U_I$.
It follows that the morphism $(\widetilde P
(\mathbf{F})_{\bf k}, u) : CY_{\bf k}^\circ \rightarrow
\AA^2_k$ is smooth on a neighborhood of $U_I$ in $CY_{\bf
k}^\circ$, so one can argue similarly as in the proof of
%%%%%%%%mm
Lemma 5.14 of \cite{ivc}.
%%%%%%%%%%%%%
\end{proof}

Fom Lemma \ref{rio} and Lemma \ref{cal},
we may rewrite
(\ref{hotrats})
\begin{equation}\label{tyo}
i^* Z^{\gamma}_{P(\mathbf{f}),X \setminus
X_0(F)}=\sum_{\sur{\delta \in F (\Gamma)}{I \in B}}
Z_{\delta, I}(T),
\end{equation}
with
\begin{equation}\label{kol}
Z_{\delta, I}(T) = ([U_I \setminus (P_{\delta}
(\mathbf{f}_I)^{-1} (0))] \Phi_{\delta, I}(T) + [\GM \times
P_{\delta} (\mathbf{f}_I)^{-1} (0)] \Psi_{\delta, I}(T),
\end{equation}
where
\begin{equation}\label{phi}
\Phi_{\delta, I}(T)= \sum_{ \sur{N_I(\mathbf{k})\in
\sigma(\delta)}{\langle N_I(\mathbf{k}), \mathbf{1}\rangle
\leq \gamma \ell_{\Gamma}(N_I (\mathbf{k}))}}
 T^{\ell_{\Gamma}(N_I (\mathbf{k}))}\, \LL^{-\sum_i \nu_i k_i},
\end{equation}
and
\begin{equation}\label{psi}
\Psi_{\delta, I}(T)= \sum_{ \sur{N_I(\mathbf{k})\in
\sigma(\delta), n >0}{\langle N_I(\mathbf{k}),
\mathbf{1}\rangle \leq \gamma \ell_{\Gamma}Ê(N_I
(\mathbf{k})) + \gamma n}}
 T^{\ell_{\Gamma}(N_I (\mathbf{k}))+ n}\, \LL^{-\sum_i \nu_i k_i}.
\end{equation}

If $\delta$ is not contained in a coordinate hyperplane,
for $\gamma$ large enough,
the inequality
\begin{equation}
\langle N_I(\mathbf{k}), \mathbf{1}\rangle 
%%%\leq \gamma \mathbf{1}\rangle ????
\leq \gamma \ell_{\Gamma}(N_I (\mathbf{k})) + \gamma n
\end{equation}
holds for every
$N_I(\mathbf{k})$ in
$\sigma(\delta)$ and every $n \geq 0$.
It follows that
\begin{equation}
\lim_{T \mapsto \infty} \Phi_{\delta, I}(T) = \lim_{T
\mapsto \infty} \Psi_{\delta, I}Ê(T)
= \chi(N_I^{-1}(\sigma(\delta))).
\end{equation}

If $\delta$ is contained in some coordinate hyperplane,
it follows from Lemma 2.10 of \cite{ivc} that
\begin{equation}
\lim_{T \mapsto \infty} \Phi_{\delta, I}(T) = \lim_{T
\mapsto \infty} \Psi_{\delta, I}(T) 
= 0.
\end{equation}

The result follows now from the
 definition of $\Psi_{P_\delta}$ and (\ref{SC}).
 \end{proof}

\begin{example} When $p = 2$ and $P = \Sigma$, one recovers
the motivic Thom-Sebastiani, cf. \cite{ts}, \cite{looi}, and \cite{barc},
in the way stated in \cite{ivc}.
When $\mathbf{f}$ is the family of coordinate
functions on the affine space $\AA^p_k$,
formula~\ref{Spsi} specializes to the one given by
Guibert~\cite{guibert}, Proposition 2.1.6.
\end{example}

\begin{remark}Restricting
to a given point $x$ of $X_0(\mathbf{f})$ and applying the
Hodge spectrum map ${\rm Sp}$ of
%%%%%%%%%%mm
 \S \kern .15em 6 of \cite{ivc}
%%%%%%%%%%%%%%
to (\ref{SfP}), one gets a formula for the
Hodge-Steenbrink spectrum (cf. \cite{St1}, \cite{St2}
\cite{Va}) of $P(\mathbf{f})$ at $x$. It is not 
immediately clear whether this formula coincides with the
one obtained by Terasoma (\cite{terasoma}, Theorem 3.6.1).
\end{remark}

\bibliographystyle{amsplain}

\end{document}